\documentclass{amsart}
\usepackage{color}

\usepackage{tipa}
\usepackage{amssymb}
\usepackage{stmaryrd}
\usepackage{graphicx}
\usepackage{amsmath,amsfonts}
\usepackage{mathtools}
\usepackage{booktabs}
\usepackage{tikz}
\usepackage{setspace}
\usepackage{float}
\usetikzlibrary{arrows,shapes,chains}
\usepackage{array}
\newcolumntype{C}[1]{>{\centering}p{#1}}

\newtheorem{theorem}{Theorem}[section]
\newtheorem{lemma}[theorem]{Lemma}
\newtheorem{proposition}[theorem]{Proposition}
\newtheorem{conjecture}{Conjecture}[section]

\theoremstyle{definition}
\newtheorem{definition}[theorem]{Definition}

\newtheorem{question}[theorem]{Question}

\numberwithin{equation}{section}

\begin{document}

%\title[short text for running head]{full title}
\title[On the Diophantine system involving pairs of triangles]{On the Diophantine system involving pairs of triangles with the same area and the same perimeter}

\author{Yangcheng Li}
 %   Address of record for the research reported here
\address{School of Mathematics and Statistics, Changsha University of Science and Technology; Hunan Provincial Key Laboratory of Mathematical Modeling and Analysis in Engineering, Changsha 410114, People's Republic of China}
\email{liyangchengmlx@163.com}

\author{Yong zhang}
%   Address of record for the research reported here
\address{School of Mathematics and Statistics, Changsha University of Science and Technology; Hunan Provincial Key Laboratory of Mathematical Modeling and Analysis in Engineering, Changsha 410114, People's Republic of China}
\email{zhangyongzju@163.com}
\thanks{This research was supported by the National Natural Science Foundation of China (Grant No. 11501052), the Natural Science Foundation of Education Department of Hunan Province (Project No. 21C0211), the Natural Science Foundation of Hunan Province (Project No. 2021JJ30699 and No. 2022JJ40464), and Hunan Provincial Key Laboratory of Mathematical Modeling and Analysis in Engineering (Changsha University of Science and Technology).}

%   General info
\date{}


\begin{abstract}
Many authors studied the problem that rational triangle pairs (triangle-parallelogram pairs) with the same area and the same perimeter. They investigated this problem by solving the rational solutions of the corresponding Diophantine equations. In this paper, we give a unified description of this problem by using the affine transformation in a rectangular coordinate system. According to the fact that two triangles with the same area and the same perimeter determine an affine transformation, this problem can be reduced to solving a specific Diophantine system. Moreover, we will give some rational solutions to this Diophantine system.
\end{abstract}

\maketitle

\section{Introduction}
A rational (Heron) triangle is a triangle with rational (integral) sides and rational (integral) area. Many authors studied two or more triangles with the same area and the same perimeter. Prielipp \cite{Prielipp} proved that there are no two such distinct right triangles, Bradley \cite{Bradley} (also see \cite{Yiu}) showed that there are no three such isosceles triangles, Hirakawa and Matsumura \cite{Hirakawa-Matsumura} obtained a unique pair of right triangle-isosceles triangle, Choudhry \cite{Choudhry2007} (also see \cite{Aassila,Bremner,Kramer-Luca,Luijk}) investigated two or more distinct Heron (rational) triangles, Juyal and Moody \cite{Juyal-Moody} studied Heron triangle and right triangle pairs. More results can be found in \cite{Lichtenberg,Wares,Choudhry-Zargar}.

For triangle and quadrilateral pairs, Guy \cite{Guy} proved that there is no non-degenerate such integral right triangle-rectangle pair. Several authors studied other cases, such as Zhang \cite{Zhang} (integral right triangle-parallelogram and Heron triangle-parallelogram pairs), Lals\'{\i}n and Ma \cite{Lalsin-Ma} ($\theta$-triangle-$\omega$-parallelogram pairs). For more results, we can refer to \cite{Bremner-Guy,Chern,Das-Juyal-Moody,Zhang-Peng,Zhang-Peng-Wang,Zhang-Zargar,Zhang-Zargar1}.

Suppose the areas of two polygons are $A_1,A_2$, and the perimeters are $P_1,P_2$, respectively. We say that the areas and perimeters of these two polygons are in fixed proportions $(\alpha,\beta)$ respectively, if they satisfy the following relationship
\begin{equation*}
A_2=\alpha A_1,\quad P_2=\beta P_1,
\end{equation*}
where $\alpha$ and $\beta$ are positive rational numbers. When $\alpha=\beta=1$, it corresponds to two polygons with the same area and the same perimeter.

Li and Zhang \cite{Li-Zhang1,Li-Zhang2,Li-Zhang3} obtained infinitely many pairs of rational triangles (pairs of cyclic quadrilaterals, pairs of triangles and some special quadrilaterals) with areas and perimeters in fixed proportions $(\alpha,\beta)$ respectively.

On an affine plane, the affine transformation is a one-to-one, onto function which preserve parallelism. In this paper, the affine transformation is confined to a rectangular coordinate system. Let $f$ be an affine transformation and $P$ be a point with coordinates $(x,y)$. Suppose that the point $P$ becomes $f(P)$ with coordinates $(x',y')$. We have the following transformation formula:
\begin{equation}
\begin{cases}
\begin{aligned}
&x'=a_{11}x+a_{12}y+b_1,\\
&y'=a_{21}x+a_{22}y+b_2.                                        \label{1}
\end{aligned}
\end{cases}
\end{equation}
Let $T$ be a triangle and $f(T)$ be the image of $T$. Let $A$ and $A'$ be the areas of $T$ and $f(T)$ respectively, then we have $A'=\vert\det(M_f)\vert A$, where $M_f$ is the transformation matrix of $f$, i.e.,
$$
M_f=\left[\begin{matrix}
a_{11} & a_{12} \\
a_{21} & a_{22}
\end{matrix}\right].$$

The structure of this paper is as follows. In Section 2, we give the definition of the weakly metric equivalence of two rational triangles and the set $VL_2(\mathbb{Q})$ to describe the affine transformation between two weakly metric equivalent rational triangles. In Section 3, the set $VL_2(\mathbb{Q})$ is used to uniformly describe the existed conclusions that two rational triangles are weakly metric equivalent, and to get more such pairs. In Section 4, the affine transformation is also used for describing the problem that triangle and parallelogram pairs having the same area and the same perimeter. In Section 5, we consider the rational triangle pairs or rational triangle-parallelogram pairs with areas and perimeters in fixed proportions $(\alpha,\beta)$, respectively, where $\alpha$ and $\beta$ are positive rational numbers.

\section{Equivalence relation and the set $VL_2(\mathbb{Q})$}
The purpose of this paper is to study the problem of two rational triangles with the same area and the same perimeter. It is useful to define the following equivalence relation.
\begin{definition}
Two rational triangles are weakly metric equivalent if they have the same area and the same perimeter. In particular, two rational triangles are metric equivalent if they are congruent.
\end{definition}
Suppose that the rational triangles $T_1$ and $T_2$ are weakly metric equivalent, we write it as $T_1\sim T_2$. If the rational triangles $T_1$ and $T_2$ are metric equivalent, we write it as $T_1=T_2$.

For a rational triangle $T$, we wish to find an affine transformation $f$ such that $f(T)$ is weakly equivalent to $T$, i.e., $f(T)\sim T$. Let $V_1,V_2$ and $V_3$ be the three vertices of the rational triangle $T$. In a rectangular coordinate system, we put the rational triangle $T$ in the first quadrant. The vertex $V_1$ is at the origin of the coordinate system, the vertex $V_2$ is on the positive $x$ axis and the vertex $V_3$ is above the $x$ axis (see Figure 1). Let $f(V_1),f(V_2)$ and $f(V_3)$ be the three vertices of the rational triangle $f(T)$. We can make $f(V_1)$ coincide with $V_1$, that is still at the origin. $f(V_2)$ is on the positive $x$ axis and $f(V_3)$ is above the $x$ axis.
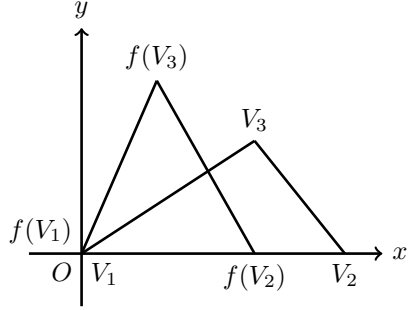
\begin{figure}[h]
\begin{tikzpicture}
\draw[line width=1pt, ->] (-0.7,0) -- (4,0) node[right] {$x$};
\draw[line width=1pt, ->] (0,-0.7) -- (0,3) node[above] {$y$};
\draw (0, 0) node[below left] {$O$};
\draw (0, 0) node[above left] {$f(V_1)$};
\draw (0, 0) node[below right] {$V_1$};
\draw[line width=1pt] (0, 0) -- (2.3, 1.5) node[above] {$V_3$};
\draw[line width=1pt] (2.3, 1.5) -- (3.5, 0) node[below] {$V_2$};
\draw[line width=1pt] (0, 0) -- (1, 2.3) node[above] {$f(V_3)$};
\draw[line width=1pt] (1, 2.3) -- (2.3, 0) node[below] {$f(V_2)$};
\end{tikzpicture}
\caption{Rational triangles $T$ and $f(T)$.}
\end{figure}

In other words, the rational triangles to be discussed are confined to the following set
\[X(T)=\{T~\vert~V_1=(0,0),V_2=(r,0),V_3=(s,t),r\in\mathbb{Q^*},s\in\mathbb{Q}_{\geq0},t\in\mathbb{Q^*}\},\]
where $\mathbb{Q^*}=\{x\in\mathbb{Q}~\vert~x>0\}$ and $\mathbb{Q}_{\geq0}=\{x\in\mathbb{Q}~\vert~x\geq0\}$. It is not difficult to find that to make $f(T)\in X(T)$, we only need $b_1=b_2=a_{21}=0$ in the transformation formula (\ref{1}). In addition, to make the areas of $T$ and $f(T)$ be the same, as long as the determinant $\det(M_f)$ of $M_f$ is equal to $\pm1$. We just take $\det(M_f)=1$. Therefore, $M_f$ can be assumed to be in the following form
$$
M_f=\left[\begin{matrix}
\frac{1}{a} & b \\
0 & a \\
\end{matrix}\right],$$
where $a\in\mathbb{Q^*},b\in\mathbb{Q}$. Thus, the transformation formula (\ref{1}) becomes
\begin{equation*}
\begin{cases}
\begin{aligned}
&x'=\frac{1}{a}x+by,\\
&y'=ay.
\end{aligned}
\end{cases}
\end{equation*}
Let
\[UL_2(\mathbb{Q})=\bigg\{\left[\begin{matrix}
\frac{1}{a} & b \\
0 & a \\
\end{matrix}\right]\bigg\vert~a\in\mathbb{Q^*},b\in\mathbb{Q}\bigg\}.\]
It is easy to check that the set $UL_2(\mathbb{Q})$ forms a group under the multiplication of matrices. Let
\[U=\{f~\vert~M_f\in UL_2(\mathbb{Q})\}.\]
Obviously, an operation of $U$ on $X(T)$ is given by the map
\[U\times X(T)\rightarrow X(T),\quad (f,T)\mapsto f(T).\]
Suppose the area of $T\in X(T)$ is $A_T$, then the orbit of $T$ in $X(T)$ is the set
\[O_T=\{T'\in X(T)~\vert~\text{the area of}~T'~\text{is}~A_T\}.\]
But the group $UL_2(\mathbb{Q})$ is not the set that we need, because we cannot guarantee that $T$ and $f(T)$ have the same perimeter.

The following definition is derived from the weakly metric equivalence of two rational triangles.
\begin{definition}
An affine transformation $f$ is called weak metric affine transformation about rational triangle $T$, if there is a rational triangle $T$ such that $f(T)$ is weakly metric equivalent to $T$. The corresponding transformation matrix $M_f$ is called weak metric transformation matrix.
\end{definition}
Let $D$ be the set of all weak metric affine transformations. Correspondingly, let $DL_2(\mathbb{Q})$ be the set of all weak metric transformation matrices. We take the intersection of the sets $UL_2(\mathbb{Q})$ and $DL_2(\mathbb{Q})$, i.e.,
\[VL_2(\mathbb{Q})=\{M_f~\vert~ M_f\in UL_2(\mathbb{Q})\cap DL_2(\mathbb{Q})\},\]
which is the set we will discuss next. Correspondingly, let
\[V=\{f~\vert ~M_f\in VL_2(\mathbb{Q})\}.\]

By Definition 2.1 and the set $V$, we have
\begin{proposition}
If $T_1\sim T_2$, then there is a $f\in V$ such that $f(T_1)=T_2$.
\end{proposition}
The set $VL_2(\mathbb{Q})$ is not empty, because $I=\left[\begin{matrix}
1 & 0 \\
0 & 1 \\
\end{matrix}\right]\in VL_2(\mathbb{Q})$. Obviously, for any rational triangle $T$, $f_I(T)=T$. Hence, the set $VL_2(\mathbb{Q})$ is a monoid.

For any $f\in V$, there is a rational triangle $T$ such that $f(T)\sim T$. Obviously, $f^{-1}(f(T))=T$ is weakly metric equivalent to $f(T)$. Therefore, for any $M_f\in VL_2(\mathbb{Q})$, $M^{-1}_f\in VL_2(\mathbb{Q})$.

However, the set $VL_2(\mathbb{Q})$ may not form a group, because it is not clear whether the composition of weak metric affine transformation is closed. Taking any $f,g\in V$ such that $f(T_1)\sim T_1$ and $g(T_2)\sim T_2$. The question is whether we can find a rational triangle $T$ such that $g(f(T))\sim T$ or $f(g(T))\sim T$. If there is a positive answer to this question, the set $VL_2(\mathbb{Q})$ can form a group.

Taking any $f\in V$, there is a rational triangle $T$ such that $f(T)\sim T$. Assuming that rational triangle $T'$ is similar to $T$, then $f(T')\sim T'$.

For the set $VL_2(\mathbb{Q})$ or $V$, we have the following two basic questions:

(1) Given a rational triangle $T$, is there a $f\in V$ such that $f(T)\sim T$?

(2) Given any $f\in V$ such that $f(T)\sim T$, is there a rational triangle $T'\in X(T)$ that is not similar to $T$ such that $f(T')\sim T'$?

To solve question (1), we give the following two lemmas:
\begin{lemma}\cite{Choudhry2007}
Given any scalene rational triangle $T_1$, whose sides do not, in any order, satisfy the Diophantine equation
\begin{align*}
&x^4-x^3y-2x^3z+3x^2yz-xy^3+3xy^2z\\
&-6xyz^2+3xz^3+y^4-2y^3z+3yz^3-2z^4=0,
\end{align*}
then there exists a second rational triangle $T_2$ such that $T_1\sim T_2$.
\end{lemma}
\begin{lemma}\cite{Choudhry2007}
Given an arbitrary positive integer $k\geq2$, there exist $k$ scalene rational triangles are weakly metric equivalent.
\end{lemma}
Lemma 2.4 gives a positive answer to question (1), and Lemma 2.5 shows that there can be an infinite number of such $f\in V$. By Lemma 2.5, we give the concept of weak similarity of weakly metric affine transformation.
\begin{definition}
Let $f,g\in V$, if there is a rational triangle $T$ and $h\in V$ such that $f(T)=h(g(T))$, then $f$ and $g$ are weakly similar, we write it as $f\sim g$.
\end{definition}
Obviously, weak similarity is an equivalence relation. All weakly metric affine transformations $f$ about the same rational triangle $T$ are weakly similar and vice versa.

For question (2), it is easy to see that $f(T)\sim T$ for any rational triangle $T$ when $f=f_I$. For a general $f\in V$, there should be $f(T)\sim T$ only for some specific rational triangles $T$. Therefore, we have the following conjecture:
\begin{conjecture}
If $f\in V, f\neq f_I$, there are only a finite number of dissimilar rational triangles $T\in X(T)$ such that $f(T)\sim T$.
\end{conjecture}

It is not clear whether the set $V$ is closed under the composition of weakly metric affine transformation. However, in some special cases, this is possible. By lemma 2.5, we can get a sequence of rational triangles $\{T_n\},n\geq1$, in which any two rational triangles are weakly metric equivalent. By Proposition 2.3, we can use $\{T_n\}$ to construct a sequence $\{f_n\}$ such that $f_i(T_i)=T_{i+1},i\geq1$. We have the following proposition:
\begin{proposition}
If $f_1,\cdots,f_k\in\{f_n\}$, $k\geq1$, then $f_k\cdots f_1\in V$. Generally, $f_k\cdots f_i\in V$, $k\geq i\geq1$.
\end{proposition}

We have a certain grasp of the properties of the set $VL_2(\mathbb{Q})$, but the conditions satisfied by $a$ and $b$ are not specific. Next, we show that the values of $a$ and $b$ are determined by a Diophantine system.

For a rational triangle $T\in X(T)$, the coordinates of the three vertices of $T$ are
\[V_1=(0,0),~~V_2=(r,0),~~V_3=(s,t),\]
where $r\in\mathbb{Q^*},s\in\mathbb{Q}_{\geq0},t\in\mathbb{Q^*}$. Thus, the sides of the rational triangle $T$ are
\[\vert V_1V_2\vert=r,~~\vert V_1V_3\vert=\sqrt{s^2+t^2},~~\vert V_2V_3\vert=\sqrt{\left(s-r\right)^2+t^2}.\]
Since the sides of rational triangle $T$ are rationals, we thus have
\[s^2+t^2=\square,~~ \left(s-r\right)^2+t^2=\square,~~r\in\mathbb{Q^*},s\in\mathbb{Q}_{\geq0},t\in\mathbb{Q^*}.\]
By the transformation $f$, the coordinates of the three vertices of $f(T)$ are
\[f(V_1)=(0,0),~~ f(V_2)=\left(\frac{r}{a},0\right),~~ f(V_3)=\left(\frac{s}{a}+bt,at\right).\]
Then, the sides of the rational triangle $f(T)$ are
\begin{align*}
&\vert f(V_1)f(V_2)\vert=\frac{r}{a},~~ \vert f(V_1)f(V_3)\vert=\sqrt{\left(\frac{s}{a}+bt\right)^2+\left(at\right)^2},\\
&\vert f(V_2)f(V_3)\vert=\sqrt{\left(\frac{s-r}{a}+bt\right)^2+(at)^2}.
\end{align*}
To ensure that the sides of $f(T)$ are rational, we have
\[\left(\frac{s}{a}+bt\right)^2+\left(at\right)^2=\square,~~\left(\frac{s-r}{a}+bt\right)^2+(at)^2=\square.\]
Since $T$ and $f(T)$ are weakly metric equivalent, they have the same perimeter, i.e.,
\[r+\sqrt{s^2+t^2}+\sqrt{\left(s-r\right)^2+t^2}=\frac{r}{a}+\sqrt{\left(\frac{s}{a}+bt\right)^2+(at)^2}+\sqrt{\left(\frac{s-r}{a}+bt\right)^2+\left(at\right)^2}.\]
Therefore, for any rational triangle $T\in X(T)$, its corresponding weakly metric affine transformation $f$ is completely characterized by the following Diophantine system
\begin{equation}
\begin{cases}
\begin{aligned}
&s^2+t^2=w_1^2,\\
&\left(s-r\right)^2+t^2=w_2^2,\\
&\left(\frac{s}{a}+bt\right)^2+\left(at\right)^2=w_3^2,\\
&\left(\frac{s-r}{a}+bt\right)^2+(at)^2=w_4^2,\\
&r+w_1+w_2=\frac{r}{a}+w_3+w_4,                                               \label{2}
\end{aligned}
\end{cases}
\end{equation}
where $b\in\mathbb{Q}$, $s\in\mathbb{Q}_{\geq0}$, $a,r,t,w_1,w_2,w_3,w_4\in\mathbb{Q^*}$.

Assume that $r,s$ and $t$ satisfy the first two equations, then $(a,b)=(1,0)$ is a solution of Eq. (\ref{2}). The corresponding affine transformation is $f_I$. Moreover, $(a,b)=\left(1,\frac{r-2s}{t}\right)$ is also a solution of Eq. (\ref{2}).

Suppose $(a,b,r,s,t)$ is any solution of Eq. (\ref{2}) and $k$ is any positive rational number, then $(a,b,kr,ks,kt)$ is also a solution of Eq. (\ref{2}), which means that Eq. (\ref{2}) is homogeneous with respect to $r,s$ and $t$.

Let $S$ be the solution set of Eq. (\ref{2}). In particular, let $S_{(r,s,t)}$ be the solution set with respect to $(a,b)$. Thus, the set $VL_2(\mathbb{Q})$ can be expressed as
\[VL_2(\mathbb{Q})=\{M~\vert~ M\in UL_2(\mathbb{Q}),(a,b)\in S_{(r,s,t)}\},\]
which gives an accurate description of the set $VL_2(\mathbb{Q})$.

Now we discuss whether the composition of weakly metric affine transformations in the set $V$ is closed according to the solution set $S$. The question can be described as: Suppose $(a_1,b_1,r_1,s_1,t_1)\in S$ and $(a_2,b_2,r_2,s_2,t_2)\in S$, if $(a,b)=\left(a_1a_2,\frac{b_2}{a_1}+a_2b_1\right)$ or $(a,b)=\left(a_1a_2,\frac{b_1}{a_2}+a_1b_2\right)$, does Eq. (\ref{2}) have a rational solution? If there is a positive answer to this question, the set $VL_2(\mathbb{Q})$ can obviously form a group. Thus, we have the following question.
\begin{question}
Does the set $VL_2(\mathbb{Q})$ form a group?
\end{question}

\section{Concrete weakly metric affine transformation}
In this section, we will discuss two problems. One is using the set $VL_2(\mathbb{Q})$ to uniformly describe the existed conclusions that two rational triangles with the same area and the same perimeter. The other is applying the set $VL_2(\mathbb{Q})$ to find more weakly metric equivalent rational triangles.

\subsection{Existed results}
We first give a unified description of the existed results, focusing on the following four cases.

\underline{Case 3.1.1.} $T$ and $f(T)$ are right triangles.

Without losing generality, we take $r=2km,s=0,t=k(m^2-1),k>0,m>1$. Then the coordinates of the three vertices of $T$ are
\[V_1=(0,0),~~V_2=(2km,0),~~V_3=(0,k(m^2-1)).\]
The coordinates of the three vertices of $f(T)$ are
\[f(V_1)=(0,0),~~f(V_2)=\left(\frac{2km}{a},0\right),~~f(V_3)=\left(bk(m^2-1),ak(m^2-1)\right).\]
Since $f(T)$ is a right triangle, we can take $b=0$, then $w_3=ak(m^2-1)$. Eq. (\ref{2}) now becomes
\begin{equation}
\begin{cases}
\begin{aligned}
&\frac{4k^2m^2}{a^2}+a^2k^2(m^2-1)^2=w_4^2,\\
&2km(m+1)=\frac{2km}{a}+ak(m^2-1)+w_4.                                               \label{3}
\end{aligned}
\end{cases}
\end{equation}
Eliminating $w_4$ in Eq.\ (\ref{3}), we have
\begin{equation*}
(a-1)((m^2-1)a-2m)=0.
\end{equation*}
Then $a=1$ or $a=\frac{2m}{m^2-1}$. Obviously, when $a=1$ or $a=\frac{2m}{m^2-1}$, $f(T)$ and $T$ are congruent. Therefore, there are no two distinct right triangles are weakly metric equivalent, which is obtained by Prielipp \cite{Prielipp}.

\underline{Case 3.1.2.} $T$ and $f(T)$ are isosceles triangles.

\underline{Case 3.1.2.1.} Take $r=4km,s=2km,t=k(m^2-1),k>0,m>1$, then the coordinates of the three vertices of $T$ are
\[V_1=(0,0),~~V_2=(4km,0),~~V_3=(2km,k(m^2-1)).\]
The coordinates of the three vertices of $f(T)$ are
\[f(V_1)=(0,0),~~f(V_2)=\left(\frac{4km}{a},0\right),~~f(V_3)=\left(\frac{2km}{a}+bk(m^2-1),ak(m^2-1)\right).\]
Since $f(T)$ is an isosceles triangle, we can take $b=0$, then $w_3=w_4$. Eq. (\ref{2}) now becomes
\begin{equation}
\begin{cases}
\begin{aligned}
&\frac{4k^2m^2}{a^2}+a^2k^2(m^2-1)^2=w_4^2,\\
&2k(m+1)^2=\frac{4km}{a}+2w_4.                                               \label{4}
\end{aligned}
\end{cases}
\end{equation}
Eliminating $w_4$ in Eq.\ (\ref{4}), we have
\begin{equation*}
(a-1)((m-1)^2a^2+(m-1)^2a-4m)=0.
\end{equation*}
Then
\[a=1~~\text{or}~~a=\frac{1-m\pm\sqrt{m^2+14m+1}}{2(m-1)}.\]
When $a=1$, $f=f_I$, this is a trivial case. Since $a\in\mathbb{Q^*}$, then $a$ has at most one nontrivial positive rational solution, which leads to the fact that there are no three isosceles triangles are weakly metric equivalent at the same time.

\underline{Case 3.1.2.2.} Take $r=2k(m^2-1),s=k(m^2-1),t=2km,k>0,m>1$, then the coordinates of the three vertices of $T$ are
\[V_1=(0,0),~~V_2=(2k(m^2-1),0),~~V_3=(k(m^2-1),2km).\]
The coordinates of the three vertices of $f(T)$ are
\[f(V_1)=(0,0),~~f(V_2)=\left(\frac{2k(m^2-1)}{a},0\right),~~f(V_3)=\left(\frac{k(m^2-1)}{a}+2bkm,2akm\right).\]
Since $f(T)$ is an isosceles triangle, we can take $b=0$, then $w_3=w_4$. Eq. (\ref{2}) now becomes
\begin{equation}
\begin{cases}
\begin{aligned}
&\frac{k^2(m^2-1)^2}{a^2}+4a^2k^2m^2=w_4^2,\\
&4km^2=\frac{2k(m^2-1)}{a}+2w_4.                                               \label{5}
\end{aligned}
\end{cases}
\end{equation}
Eliminating $w_4$ in Eq.\ (\ref{5}), we have
\begin{equation*}
(a-1)(a^2+a-(m^2-1))=0.
\end{equation*}
Then
\[a=1~~\text{or}~~a=\frac{-1\pm\sqrt{4m^2-3}}{2}.\]
Hence, $a$ has at most one nontrivial positive rational solution, then there are no three isosceles triangles are weakly metric equivalent at the same time, this case was studied by Bradley \cite{Bradley} (also see \cite{Yiu}).

\underline{Case 3.1.3.} $T$ is a right triangle and $f(T)$ is an isosceles triangle.

Take $r=2km,s=0,t=k(m^2-1),k>0,m>1$. Then the coordinates of the three vertices of $T$ are
\[V_1=(0,0),~~V_2=(2km,0),~~V_3=(0,k(m^2-1)).\]
The coordinates of the three vertices of $f(T)$ are
\[f(V_1)=(0,0),~~f(V_2)=\left(\frac{2km}{a},0\right),~~f(V_3)=\left(bk(m^2-1),ak(m^2-1)\right).\]
Since $f(T)$ is an isosceles triangle, we can take $bk(m^2-1)=\frac{km}{a}$, then $w_3=w_4$. Eq. (\ref{2}) now becomes
\begin{equation}
\begin{cases}
\begin{aligned}
&\frac{k^2m^2}{a^2}+a^2k^2(m^2-1)^2=w_4^2,\\
&2km(m+1)=\frac{2km}{a}+2w_4.                                               \label{6}
\end{aligned}
\end{cases}
\end{equation}
Eliminating $w_4$ in Eq.\ (\ref{6}), we have
\begin{equation}
(m+1)(m-1)^2a^3-m^2(m+1)a+2m^2=0.                                               \label{7}
\end{equation}
The discriminant of Eq.\ (\ref{7}) with respect to $a$ is
\[4(m^2-1)^2m^4(m^4+2m^3-26m^2+54m-27)>0.\]
Eq.\ (\ref{7}) is equivalent to
\begin{equation*}
(m+1)a((m-1)^2a^2-m^2)+2m^2=0,
\end{equation*}
then $a<\frac{m}{m-1}$.

It is easy to verify that $(a,m)=\left(\frac{8}{3},\frac{16}{11}\right)$ $(k=121)$ is the solution of Eq.\ (\ref{7}), which corresponds to the unique pair of right triangle and isosceles triangle obtained by Hirakawa and Matsumura \cite{Hirakawa-Matsumura}, i.e., $(135,352,377), (132,366,366)$. Unfortunately, we cannot prove that this solution is the only nontrivial positive rational solution of Eq.\ (\ref{7}).

\underline{Case 3.1.4.} $T$ and $f(T)$ are rational triangles.

In this case, we only give the affine transformation corresponding to a parameter solution obtained in \cite{Choudhry2007}. The sides of $T$ and $f(T)$ are given respectively by
\begin{align*}
a_1=&~(R^2+1)(S^2+S+1),\\
b_1=&~S(S+1)(R^2+S^2+S+1),\\
c_1=&~R^2+(S^2+S+1)^2,
\end{align*}
and
\begin{align*}
a_2=&~(R^2+S^2)(S^2+S+1),\\
b_2=&~(S+1)(R^2+S^2+S+1),\\
c_2=&~R^2S^2+(S^2+S+1)^2,
\end{align*}
where $R,S$ are arbitrary positive rational numbers. Take
\begin{align*}
r=&~R^2+(S^2+S+1)^2,\\
s=&~\frac{S(S+1)(R^2+S^2+S+1)(S^2+R+S+1)(S^2-R+S+1)}{S^4+2S^3+R^2+3S^2+2S+1},\\
t=&~\frac{(2S^2+2S+2)(R^2+S^2+S+1)RS(S+1)}{S^4+2S^3+R^2+3S^2+2S+1}.
\end{align*}
Then
\begin{align*}
a=&~\frac{S^4+2S^3+R^2+3S^2+2S+1}{R^2S^2+S^4+2S^3+3S^2+2S+1},\\
b=&~\frac{(S-1)b(R,S)}{2SR(S^4+2S^3+R^2+3S^2+2S+1)(R^2S^2+S^4+2S^3+3S^2+2S+1)},
\end{align*}
where
\begin{align*}
b(R,S)=&~R^6S^2-(S^2+S+1)(S^4+S^3-S^2+S+1)R^4\\
&-(2S^2+3S+2)(S^2+S+1)^3R^2-(S^2+S+1)^5.
\end{align*}

\subsection{Other solutions of Eq. (\ref{2})}
We can get more solutions by directly solving Eq. (\ref{2}). We are concerned with the following two cases.

\underline{Case 3.2.1.} $s=0$, $T$ is a right triangle.

Take $r=2km,s=0,t=k(m^2-1),k>0,m>1$, then
\[w_1=k(m^2-1),~~w_2=k(m^2+1).\]
From the third equation of Eq. (\ref{2}), we have
\[k^2(m^2-1)^2(a^2+b^2)=w_3^2,\]
then
\[a=2cp,~~b=c(p^2-1),~~w_3=ck(m^2-1)(p^2+1),~~c>0,~~p>1.\]
Eq. (\ref{2}) now becomes
\begin{equation}
\begin{cases}
\begin{aligned}
&\left(ck(p^2-1)(m^2-1)-\frac{km}{cp}\right)^2+4c^2p^2k^2(m^2-1)^2=w_4^2,\\
&2km(m+1)=\frac{km}{cp}+ck(m^2-1)(p^2+1)+w_4.                                               \label{8}
\end{aligned}
\end{cases}
\end{equation}
Eliminating $w_4$ in Eq.\ (\ref{8}), we have
\begin{equation}
p(p^2+1)(m^2-1)c^2-p(m^2+(p+1)m-p)c+m=0,                                               \label{9}
\end{equation}
then
\[c=\frac{p((m-1)p+m(m+1))\pm\sqrt{\Delta(c)}}{2p(p^2+1)(m^2-1)},\]
where
\[\Delta(c)=p^2m^4-2p(p^2-p+2)m^3+p^2(p^2+1)m^2-2p(p^3-p^2-2)m+p^4.\]
To get the rational solutions of Eq.\ (\ref{9}) with respect to $c$, it needs $\Delta(c)$ be a perfect square, say $w^2$. Using a method described by Fermat (see \cite[p.~639]{Dickson}), it's easy to get the following solutions
\begin{equation*}
\begin{split}
m=\frac{p^2+p+1}{p^2},~~m=\frac{p^2-p+1}{p(p-2)},~~m=\frac{p^3-p^2-1}{p^2(p+1)}
\end{split}
\end{equation*}
such that $\Delta(c)$ is a perfect square. Take $m=\frac{p^2+p+1}{p^2}$, by Eq.\ (\ref{9}), we get
\[c=\frac{p}{p^2+1},~~c=\frac{p^2+p+1}{2p^3+3p^2+2p+1}.\]
Take $c=\frac{p}{p^2+1}$, then
\begin{align*}
&s=0,~~r=\frac{2k(p^2+p+1)}{p^2},~~t=\frac{k(p+1)(2p^2+p+1)}{p^4},\\
&a=\frac{2p^2}{p^2+1},~~b=\frac{p(p^2-1)}{p^2+1},\\
&w_1=\frac{k(p+1)(2p^2+p+1)}{p^4},~~w_2=\frac{k(p^2+1)(2p^2+2p+1)}{p^4},\\
&w_3=\frac{k(p+1)(2p^2+p+1)}{p^3},~~w_4=\frac{(p^4+2p^3+4p^2+2p+1)k}{p^4},
\end{align*}
where $k>0,p>1$, this is a parametric solution of Eq. (\ref{2}).

Let us consider the rational points on the quartic curve
\[\mathcal{C}: w^2=\Delta(c).\]
By the map $\varphi$:
\begin{align*}
X=&~\frac{3p(p(p^2+1)m^2-6(p^3-p^2-2)m+6p^3+6pw)}{m^2},\\
Y=&-\frac{Y(m,w)}{m^3},
\end{align*}
where
\begin{align*}
Y(m,w)=&~54p(p^2((p^2-p+2)m^3-p(p^2+1)m^2+3(p^3-p^2-2)m-2p^3)\\
&+((p^3-p^2-2)m-2p^3)w).
\end{align*}
We can transform $\mathcal{C}$ into the elliptic curve
\begin{align*}
\mathcal{E}_{p}:~Y^2=~&X^3-27p^2(p^6-12p^5+38p^4-36p^3+49p^2-24p+48)X\\
&+54p^4(p^3-6p^2+p-12)(p^5-12p^4+38p^3-36p^2+p-24).
\end{align*}
The discriminant of $\mathcal{E}_{p}$ is
\[\Delta(p)=544195584p^6(p^4-10p^3+17p^2+8p+16)(p-1)^2(p^2+1)^2.\]
Hence, $\Delta(p)\neq0$, so $\mathcal{E}_{p}$ is non-singular.

It is easy to check that the elliptic curve $\mathcal{E}_{p}$ contains two rational points
\begin{align*}
Q=(3p(p^3-6p^2+p-12),0),~~P=(3p^2(p^2-6p+1),108p^2(p+1)).
\end{align*}
By the inverse map $\varphi^{-1}$, we have
\begin{align*}
m=\frac{6p(Yp-3(p^3-p^2-2)X+9p^2(p^5-7p^4+7p^3-15p^2-2))}{(X-3p^2(p^2-6p+1))(X-3p^2(p^2+6p+1))}.
\end{align*}
When $p=2$, we get
\begin{align*}
\mathcal{E}_{2}:~Y^2=~X^3-21168X+494208.
\end{align*}
By the package of Magma, the rank of $\mathcal{E}_{2}$ is $1$. From the theorem of Silverman (see p. 457 of \cite{Silverman}) about the rank of elliptic curve: If an elliptic curve $E$ defined over $\mathbb{Q}(t)$ for some values of $t$ has positive rank, for all but finitely many $t\in\mathbb{Q}$ it has positive rank. Hence, the elliptic curve $\mathcal{E}_{p}$, for all but finitely many $p>1$, has positive rank. Therefore, we can get more parametric solutions by the group law, such as
\begin{align*}
[2]P=\bigg(&\frac{3(p^6-4p^5+38p^4-16p^3+25p^2-12p+12)}{(p+1)^2},\\
&\frac{108(2p^2+1)(2p^5-4p^4+5p^3-2p^2+3p-2)}{(p+1)^3}\bigg).
\end{align*}
From the rational point $[2]P$, we obtain
\[m=\frac{p(p+1)}{p^2-p+1},\]
then
\[c=\frac{p+1}{2p^2+1},~~c=\frac{p^2-p+1}{2p^3-p^2+2p-1}.\]
Take $c=\frac{p+1}{2p^2+1}$, then
\begin{align*}
&s=0,~~r=\frac{2kp(p+1)}{p^2-p+1},~~t=\frac{k(2p-1)(2p^2+1)}{(p^2-p+1)^2},\\
&a=\frac{2(p+1)p}{2p^2+1},~~b=\frac{(p+1)^2(p-1)}{2p^2+1},\\
&w_1=\frac{k(2p-1)(2p^2+1)}{(p^2-p+1)^2},~~w_2=\frac{k(2p^4+4p^2-2p+1)}{(p^2-p+1)^2},\\
&w_3=\frac{(p^2+1)(2p-1)k(p+1)}{(p^2-p+1)^2},~~w_4=\frac{kp(5p^2-2p+2)}{(p^2-p+1)^2},
\end{align*}
where $k>0,p>1$.

\underline{Case 3.2.2.} $s\neq0$, $T$ is a rational triangle.

Take $s=2km,t=k(m^2-1),k>0,m>1$, then
\[w_1=k(m^2+1).\]
From the second equation of Eq. (\ref{2}), we have
\begin{align*}
r=\frac{2k((m^2+1)c-2m)}{c^2-1},~~w_2=\frac{k((m^2+1)c^2-4mc+m^2+1)}{c^2-1},
\end{align*}
where $c$ is a rational number.

From the third equation of Eq. (\ref{2}), we have
\[\left(\frac{2km}{a}+bk(m^2-1)\right)^2+a^2k^2(m^2-1)^2=w_3^2.\]
Let
\begin{align*}
\frac{2km}{a}+bk(m^2-1)=d(p^2-1),~~ak(m^2-1)=2dp,
\end{align*}
then
\begin{align*}
a=\frac{2dp}{k(m^2-1)},~~b=\frac{p(p^2-1)d^2-k^2m(m^2-1)}{dkp(m^2-1)},~~w_3=d(p^2+1).
\end{align*}
Eq. (\ref{2}) now becomes
\begin{equation}
\begin{cases}
\begin{aligned}
&\left(\frac{d^2p(p^2-1)(c^2-1)-(m^2-1)((m^2+1)c-2m)k^2}{(c^2-1)dp}\right)^2+4d^2p^2=w_4^2,\\
&\frac{2k((m^2+1)c-2m)}{c-1}=\frac{(m^2-1)((m^2+1)c-2m)k^2}{(c^2-1)dp}+d(p^2+1)+w_4                               \label{2-1}
\end{aligned}
\end{cases}
\end{equation}
Eliminating $w_4$ in Eq.\ (\ref{2-1}), we have
\begin{equation}
\begin{split}
&dp((p^2+1)d-k(m^2+1))c^2\\
&-k(m-1)(p((m+1)p+m-1)d-k(m+1)(m^2+1))c\\
&-p(p^2+1)d^2+kp((m^2-1)p+2m)d-2k^2m(m^2-1)=0,                                \label{2-2}
\end{split}
\end{equation}
then
\[c=\frac{(m^4-1)k^2-p((m^2-1)p+(m-1)^2)k\pm\sqrt{\Delta(c)}}{2p(k(m^2+1)-(p^2+1))},\]
where
\begin{align*}
\Delta(c)=&~(m^4-1)^2k^4-2dp(m^4-1)((m^2-1)p+(m+1)^2)k^3\\
&+d^2p((m^2-1)^2p^3+2(m^2-1)(3m^2+2m+3)p^2\\
&+(m+1)^4p+8m(m^2-1))k^2\\
&-4d^3p^2(p^2+1)((m^2-1)p+(m+1)^2)k+4d^4p^2(p^2+1)^2.
\end{align*}
To get the rational solutions of Eq.\ (\ref{2-2}) with respect to $c$, it needs $\Delta(c)$ be a perfect square, say $w^2$. Since the discriminant $\Delta(c)$ is homogeneous with respect to $k$ and $d$, we can let $d=1$. Using a method described by Fermat (see \cite[p.~639]{Dickson}), it's easy to get the following solution
\begin{equation*}
\begin{split}
k=\frac{2(m^2+1)p^2+(m+1)^2}{(m^2+1)((m^2-1)p+(m+1)^2)}
\end{split}
\end{equation*}
such that $\Delta(c)$ is a perfect square. By Eq.\ (\ref{2-2}), we get
\[c=\frac{(m+1)(m^2+1)p^2+(m-1)(m^2+1)p+2m(m+1)}{(m^2+1)((m+1)p^2-(m-1)p+m+1)}.\]
By some calculations, we obtain
\begin{align*}
&r=\frac{2((m+1)p^2-(m-1)p+m+1)(m^2+1)p}{(m+1)(2(m^2+1)p-m^2+1)},\\
&s=\frac{2(2(m^2+1)p^2+(m+1)^2)m}{(m^2+1)((m^2-1)p+(m+1)^2)},\\
&t=\frac{(m-1)(2(m^2+1)p^2+(m+1)^2)}{((m-1)p+m+1)(m^2+1)},\\
&a=\frac{2((m-1)p+m+1)(m^2+1)p}{(m-1)(2(m^2+1)p^2+(m+1)^2)},\\
&b=\frac{b(m,p)}{(m^4-1)(2(m^2+1)p^2+(m+1)^2)((m^2-1)p+(m+1)^2)p},\\
\end{align*}
where
\begin{align*}
b(m,p)=&~(m^4-1)^2p^5+2(m^2-1)(m^2+1)^3p^4+4(m+1)^2m(m^2+1)^2p^3\\
&-2(m+1)^4(m^4-1)p^2-(m+1)^4(m^2+1)^2p-(m+1)^5m(m-1),
\end{align*}
and $m>1,p>1$, this is a parametric solution of Eq. (\ref{2}).

Let us consider the rational points on the quartic curve
\[\mathcal{C}: w^2=\Delta(c).\]
By the map $\varphi$:
\begin{align*}
X=&\frac{3pX(k,w)}{k^2},\\
Y=&-\frac{108p^2(p^2+1)Y(k,w)}{k^3},
\end{align*}
where
\begin{align*}
X(k,w)=&~((m^2-1)^2p^3+2(m^2-1)(3m^2+2m+3)p^2+(m+1)^4p\\
&+8m(m^2-1))k^2-12p(p^2+1)((m^2-1)p+(m+1)^2)k\\
&+24p(p^2+1)^2-12(p^2+1)w,\\
Y(k,w)=&~(m^4-1)((m^2-1)p+(m+1)^2)k^3-((m^2-1)^2p^3\\
&+2(m^2-1)(3m^2+2m+3)p^2+(m+1)^4p+8m(m^2-1))k^2\\
&+6p(p^2+1)((m^2-1)p+(m+1)^2)k-8p(p^2+1)^2\\
&-(((m^2-1)p+(m+1)^2)k-4(p^2+1))w.
\end{align*}
We can transform $\mathcal{C}$ into the elliptic curve
\begin{align*}
\mathcal{E}_{(m,p)}:~Y^2=~X^3-27p^2A_4(m,p)X+54p^3A_6(m,p),
\end{align*}
where
\begin{align*}
A_4(m,p)=&~(m^2-1)^4p^6-4(3m^2-2m+3)(m^2-1)^3p^5\\
&+2(19m^4-20m^3+50m^2-20m+19)(m^2-1)^2p^4\\
&-4(m+1)^3(m-1)(3m^2-2m+3)^2p^3\\
&+(m+1)^2(49m^6-90m^5+223m^4-300m^3+223m^2-90m+49)p^2\\
&-8(m+1)^5(m-1)(3m^2-2m+3)p+16(m^2+3)(3m^2+1)(m^2-1)^2,\\
A_6(m,p)=&~((m^2-1)^2p^3-2(m^2-1)(3m^2-2m+3)p^2+(m+1)^4p\\
&-4(m^2-1)(3m^2-2m+3))((m^2-1)^4p^6\\
&-4(3m^2-2m+3)(m^2-1)^3p^5\\
&+2(19m^4-20m^3+50m^2-20m+19)(m^2-1)^2p^4\\
&-4(m+1)^3(m-1)(3m^2-2m+3)^2p^3\\
&+(m+1)^2(m^6+102m^5-113m^4+84m^3-113m^2+102m+1)p^2\\
&-8(m+1)^5(m-1)(3m^2-2m+3)p+32m(3m^2+2m+3)(m^2-1)^2).
\end{align*}
By some calculations, we can check that when $m\neq\frac{p+1}{p-1}$, the discriminant $\Delta(m,p)$ of $\mathcal{E}_{(m,p)}$ is not $0$, then $\mathcal{E}_{(m,p)}$ is non-singular.

It is easy to check that the elliptic curve $\mathcal{E}_{(m,p)}$ contains two rational points
\begin{align*}
&Q=(X_1(m,p),0),\\
&P=(X_2(m,p),216p^2(p^2+1)(m^4-1)((m^2-1)p+(m+1)^2)),
\end{align*}
where
\begin{align*}
X_1(m,p)=&~3p((m^2-1)^2p^3-2(m^2-1)(3m^2-2m+3)p^2\\
&+(m+1)^4p-4(m^2-1)(3m^2-2m+3)),\\
X_2(m,p)=&~3p((m^2-1)^2p^3+2(m^2-1)(9m^2+2m+9)p^2\\
&+(m+1)^4p+4(m^2-1)(3m^2+2m+3)).
\end{align*}
By the inverse map $\varphi^{-1}$, we have
\begin{align*}
k=-\frac{12p(p^2+1)k(X,Y)}{(X-X_1(m,p))(X-X_2(m,p))},
\end{align*}
where
\begin{align*}
k(X,Y)=&~Y+3p((m^2-1)p+(m+1)^2)(X-3p((m^2-1)^2p^3\\
&-2(m^2-1)(3m^2-2m+3)p^2+(m+1)^4p\\
&-4(m^2-1)(3m^2-2m+3))).
\end{align*}

Similar to Case 3.2.1, we can get more parametric solutions by the group law. From the rational point $[2]P$, we get
\begin{align*}
k=\frac{2(p^2+1)p((m-1)p+m+1)}{(m-1)(2(m^2+1)p^2+(m+1)^2)},
\end{align*}
then
\[c=\frac{2(m^2-1)p^2+4pm+m^2-1}{2(m^2+1)p-m^2+1}.\]
By some calculations, we obtain
\begin{align*}
&r=\frac{(p^2+1)(2(m^2+1)p-m^2+1)}{(m^2-1)p^2-(m-1)^2p+m^2-1},\\
&s=\frac{4(p^2+1)p((m-1)p+m+1)m}{(m-1)(2(m^2+1)p^2+(m+1)^2)},\\
&t=\frac{2(p^2+1)p((m^2-1)p+(m+1)^2)}{2(m^2+1)p^2+(m+1)^2},\\
&a=\frac{2(m^2+1)p^2+(m+1)^2}{(p^2+1)((m^2-1)p+(m+1)^2)},\\
&b=-\frac{b(m,p)}{2(m^2-1)p((m^2-1)p+(m+1)^2)(p^2+1)(2(m^2+1)p^2+(m+1)^2)},
\end{align*}
where
\begin{align*}
b(m,p)=&~4m(m^2-1)^2p^7-4(m^2-1)(m^4-2m^3-2m^2-2m+1)p^6\\
&+4(m+1)^2m(3m^2-2m+3)p^5+8m(m^2-1)(m^2+4m+1)p^4\\
&+4(m+1)^2m(3m^2+2m+3)p^3+3(m-1)(m+1)^5p^2\\
&+4(m+1)^4mp+(m-1)(m+1)^5,
\end{align*}
and $m>1,p>1$, this is a parametric solution of Eq. (\ref{2}).

\section{triangle and parallelogram pairs}
In this section, we consider the problem that triangle and parallelogram pairs have the same area and the same perimeter.

\subsection{The corresponding Diophantine system}
A triangle cannot be transformed into a parallelogram by affine transformation, but the parallelogram is composed of two congruent triangles (see Figure 2), which allows us to study this problem by the affine transformation.
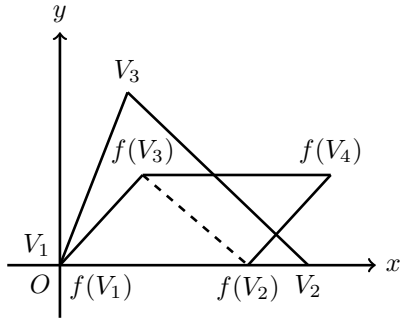
\begin{figure}[h]
\begin{tikzpicture}
\draw[line width=1pt, ->] (-0.7,0) -- (4.2,0) node[right] {$x$};
\draw[line width=1pt, ->] (0,-0.7) -- (0,3.1) node[above] {$y$};
\draw (0, 0) node[below left] {$O$};
\draw (0, 0) node[above left] {$V_1$};
\draw (0, 0) node[below right] {$f(V_1)$};

\draw[line width=1pt] (0, 0) -- (1.1, 1.2) node[above] {$f(V_3)$};
\draw[line width=1pt] (1.1, 1.2) -- (3.6, 1.2) node[above] {$f(V_4)$};
\draw[line width=1pt] (3.6, 1.2) -- (2.5, 0) node[below] {$f(V_2)$};
\draw[line width=1pt, dashed] (1.1, 1.2) -- (2.5, 0);
\draw[line width=1pt] (0, 0) -- (0.9, 2.3) node[above] {$V_3$};
\draw[line width=1pt] (0.9, 2.3) -- (3.3, 0) node[below] {$V_2$};
\end{tikzpicture}
\caption{Rational triangle $T$ and parallelogram $L$.}
\end{figure}

Let $V_1,V_2$ and $V_3$ be the three vertices of rational triangle $T$, and $T\in X(T)$. Let $f(V_1),f(V_2),f(V_3)$ and $f(V_4)$ be the four vertices of parallelogram $L$. Let $f(T)$ be a triangle with $f(V_1),f(V_2)$ and $f(V_3)$ as vertices. Let $A$ and $A'$ be the areas of rational triangle $T$ and $f(T)$, respectively. Considering the affine transformation $f$ from $T$ to $f(T)$ such that $A=2A'$, then the area of parallelogram $L$ is equal to the area of rational triangle $T$. The transformation matrix $M_f$ can be taken as
$$
M_f=\left[\begin{matrix}
\frac{1}{2a} & b \\
0 & a \\
\end{matrix}\right],$$
where $a\in\mathbb{Q^*},b\in\mathbb{Q}$. In order to ensure that rational triangle $T$ and parallelogram $L$ have the same perimeter, we have the following Diophantine system.
\begin{equation}
\begin{cases}
\begin{aligned}
&s^2+t^2=w_1^2,\\
&\left(s-r\right)^2+t^2=w_2^2,\\
&\left(\frac{s}{2a}+bt\right)^2+\left(at\right)^2=w_3^2,\\
&r+w_1+w_2=2\left(\frac{r}{2a}+w_3\right),                                               \label{4-1}
\end{aligned}
\end{cases}
\end{equation}
where $b\in\mathbb{Q}$, $s\in\mathbb{Q}_{\geq0}$, $a,r,t,w_1,w_2,w_3\in\mathbb{Q^*}$.

\subsection{Existed results}
We use Eq. (\ref{4-1}) to discuss the existed results, focusing on the following two cases.

\underline{Case 4.2.1.} $T$ is a right triangle and $L$ is a rectangle.

Take $r=2km,s=0,t=k(m^2-1),k>0,m>1$, then
\[w_1=k(m^2-1),~~w_2=k(m^2+1).\]
Since $L$ is a rectangle, we can take $b=0$. Eq. (\ref{4-1}) now becomes
\begin{equation}
2km(m+1)=2\left(\frac{km}{a}+ak(m^2-1)\right).                                              \label{4-2}
\end{equation}
Solving Eq.\ (\ref{4-2}), we get
\begin{equation*}
a=\frac{m^2+m\pm\sqrt{\Delta(a)}}{2(m^2-1)},
\end{equation*}
where
\[\Delta(a)=m^4-2m^3+m^2+4m.\]
To get the rational solutions of Eq.\ (\ref{4-2}), we need to make the discriminant $\Delta(a)$ be a perfect square, say $w^2$. By the map
\begin{align*}
X=&\frac{3(m^2+10m-12w+13)}{(m-1)^2},\\
Y=&-\frac{108(m^3-2m^2+7m+2-(m+3)w)}{(m-1)^3},
\end{align*}
We get an elliptic curve
\begin{align*}
\mathcal{E}:~Y^2=~X^3-675X+13662.
\end{align*}
The discriminant of $\mathcal{E}$ is $-3809369088$, so $\mathcal{E}$ is nonsingular.

The inverse map is
\[\begin{split}
m=~&\frac{X^2+30X-99-12Y}{(X+33)(X-39)},\\
w=~&-\frac{2(X^3-27X^2+2727X-31293-6(X+33)Y)}{(X+33)(X-39)^2}.
\end{split}\]
By the package of Magma, the rank of $\mathcal{E}$ is $0$, and $\mathcal{E}$ only has degenerate solutions
\[(X,Y)=(-33,0),(3,\pm108),(39,\pm216).\]
Therefore, there is no right triangle and rectangle pair with the same area and the same perimeter, this case was obtained by Guy \cite{Guy}.

\underline{Case 4.2.2.} $T$ is an isosceles triangle and $L$ is a rhombus.

\underline{Case 4.2.2.1.} Take $r=4km,s=2km,t=k(m^2-1),k>0,m>1$, then
\[w_1=w_2=k(m^2+1).\]
Since $L$ is a rhombus, then Eq. (\ref{4-1}) becomes
\begin{equation}
\begin{cases}
\begin{aligned}
&\left(\frac{2km}{a}\right)^2=\left(\frac{km}{a}+bk(m^2-1)\right)^2+a^2k^2(m^2-1)^2,\\
&2k(m+1)^2=\frac{8km}{a}.                                                                      \label{4-3}
\end{aligned}
\end{cases}
\end{equation}
Eliminating $a$ in Eq.\ (\ref{4-3}), we have
\begin{equation}
\begin{split}
&16(m-1)^2(m+1)^4b^2+8(m-1)(m+1)^5b\\
&-3m^6-18m^5+211m^4-572m^3+211m^2-18m-3=0.                                \label{4-4}
\end{split}
\end{equation}
The discriminant of Eq. (\ref{4-4}) with respect to $b$ is
\[\Delta(b)=(m^3+11m^2-5m+1)(m^3-5m^2+11m+1).\]
Obviously, we need to make $\Delta(b)$ be a perfect square, say $w^2$. We consider the curve
\[\mathcal{C}_1: w^2=(m^3+11m^2-5m+1)(m^3-5m^2+11m+1).\]
By the package of Magma, the Mordell-Weil rank of the Jacobian variety of $\mathcal{C}_1$ over $\mathbb{Q}$ is at most $1$, and Magma's Chabauty routines determine the only finite rational points in $\mathcal{C}_1$ are
\[(m,w)=(0,\pm1),(1,\pm8).\]
Obviously, all points are degenerate points.

\underline{Case 4.2.2.2.} Take $r=2k(m^2-1),s=k(m^2-1),t=2km,k>0,m>1$, then
\[w_1=w_2=k(m^2+1).\]
Since $L$ is a rhombus, then Eq. (\ref{4-1}) becomes
\begin{equation}
\begin{cases}
\begin{aligned}
&\left(\frac{k(m^2-1)}{a}\right)^2=\left(\frac{k(m^2-1)}{2a}+2bkm\right)^2+(2akm)^2,\\
&4km^2=\frac{4k(m^2-1)}{a}.                                                                      \label{4-5}
\end{aligned}
\end{cases}
\end{equation}
Eliminating $a$ in Eq.\ (\ref{4-5}), we have
\begin{equation}
\begin{split}
&16b^2m^4+8bm^5-3m^6+16m^4-32m^2+16=0.                                \label{4-6}
\end{split}
\end{equation}
The discriminant of Eq. (\ref{4-6}) with respect to $b$ is
\[\Delta(b)=(m^3-2m^2+2)(m^3+2m^2-2).\]
We consider the curve
\[\mathcal{C}_2: w^2=(m^3-2m^2+2)(m^3+2m^2-2).\]
The Mordell-Weil rank of the Jacobian variety of $\mathcal{C}_2$ over $\mathbb{Q}$ is at most $1$, and Magma's Chabauty routines determine the only finite rational
points in $\mathcal{C}_2$ are
\[(m,w)=(\pm1,\pm1).\]
Therefore, there is no isosceles triangle and rhombus pair with the same area and the same perimeter, this case was obtained by Zhang and Peng \cite{Zhang-Peng}.

\subsection{Other solutions of Eq. (\ref{4-1})}
We can get more solutions by directly solving Eq. (\ref{4-1}). We are concerned with the following two cases.

\underline{Case 4.3.1.} $s=0$, $T$ is a right triangle.

Take $r=2km,s=0,t=k(m^2-1),k>0,m>1$, then
\[w_1=k(m^2-1),~~w_2=k(m^2+1).\]
From the third equation of Eq. (\ref{4-1}), we have
\[k^2(m^2-1)^2(a^2+b^2)=w_3^2,\]
then
\[a=2cp,~~b=c(p^2-1),~~w_3=ck(m^2-1)(p^2+1),~~c>0,~~p>1.\]
Eq. (\ref{4-1}) now becomes
\begin{equation}
2p(p^2+1)(m^2-1)c^2-2m(m+1)pc+m=0,                                               \label{4-7}
\end{equation}
then
\[c=\frac{m(m+1)p\pm\sqrt{\Delta(c)}}{2p(p^2+1)(m^2-1)},\]
where
\[\Delta(c)=p^2m^4-2p(p^2-p+1)m^3+p^2m^2+2p(p^2+1)m.\]
To get the rational solutions of Eq.\ (\ref{4-7}) with respect to $c$, it needs $\Delta(c)$ be a perfect square, say $w^2$. Using a method described by Fermat (see \cite[p.~639]{Dickson}), it's easy to get the following solutions
\begin{equation*}
\begin{split}
m=-\frac{(p-1)^4}{4p(p^2-3p+1)}
\end{split}
\end{equation*}
such that $\Delta(c)$ is a perfect square. By Eq.\ (\ref{4-7}), we get
\[c=-\frac{(p-1)^2}{(p^2+1)(p^2-4p+1)}.\]
By some calculations, we obtain
\begin{align*}
&s=0,~~r=-\frac{k(p-1)^4}{2p(p^2-3p+1)},~~t=\frac{k(p^2+2p-1)(p^2-2p-1)(p^2-4p+1)^2}{16p^2(p^2-3p+1)^2},\\
&a=-\frac{2(p-1)^2p}{(p^2+1)(p^2-4p+1)},~~b=-\frac{(p-1)^3(p+1)}{(p^2+1)(p^2-4p+1)},
\end{align*}
where $k>0,1+\sqrt{2}<p<\frac{3+\sqrt{5}}{2}$, this is a parametric solution of Eq. (\ref{4-1}).

Let us consider the rational points on the quartic curve
\[\mathcal{C}: w^2=\Delta(c).\]
By the map $\varphi$:
\begin{align*}
X=&-\frac{3p(6(m^2-1)p^2-p(13m^2+10m+1)+6(m^2-1)+12w)}{(m-1)^2},\\
Y=&\frac{108pY(m,w)}{(m-1)^3},
\end{align*}
where
\begin{align*}
Y(m,w)=&~p((m-1)(m^2+4m+1)p^2-m(m+1)(3m+1)p\\
&+(m-1)(m^2+4m+1))-((m-1)p^2-(3m+1)p+m-1)w.
\end{align*}
We can transform $\mathcal{C}$ into the elliptic curve
\begin{align*}
\mathcal{E}_{p}:~Y^2=~&X^3-27p^2(12p^4-12p^3+25p^2-12p+12)X\\
&+54p^4(12p^2-p+12)(6p^2-p+6).
\end{align*}
The discriminant of $\mathcal{E}_{p}$ is
\[\Delta(p)=34012224p^6(4p^4-12p^3+9p^2-12p+4)(p^2+1)^4.\]
Hence, $\Delta(p)\neq0$, so $\mathcal{E}_{p}$ is non-singular.

It is easy to check that the elliptic curve $\mathcal{E}_{p}$ contains two rational points
\begin{align*}
Q=(-3p(6p^2-p+6),0),~~P=(-3p(6p^2-25p+6),-216p^2(p^2-3p+1)).
\end{align*}
By the inverse map $\varphi^{-1}$, we have
\begin{align*}
m=\frac{(X+3p(6p^2-p+6))(X-3p(6p^2-11p+6))-12pY}{(X+3p(6p^2-p+6))(X+3p(6p^2-25p+6))}.
\end{align*}
When $p=2$, we get
\begin{align*}
\mathcal{E}_{2}:~Y^2=~X^3-19872X+1403136.
\end{align*}
The rank of $\mathcal{E}_{2}$ is $1$. From the theorem of Silverman (see p. 457 of \cite{Silverman}), the elliptic curve $\mathcal{E}_{p}$, for all but finitely many $p>1$, has positive rank. Therefore, we can get more parametric solutions by the group law. From the rational point $[2]P$, we obtain
\[m=\frac{(p^2-4p+1)^2}{(p^2-2p-1)(p^2+2p-1)},\]
then
\[c=-\frac{p^2-4p+1}{4(p-1)^2p}.\]
By some calculations, we have
\begin{align*}
&s=0,~~r=\frac{2k(p^2-4p+1)^2}{(p^2+2p-1)(p^2-2p-1)},~~t=-\frac{16k(p-1)^4p(p^2-3p+1)}{(p^2+2p-1)^2(p^2-2p-1)^2},\\
&a=-\frac{p^2-4p+1}{2(p-1)^2},~~b=-\frac{(p+1)(p^2-4p+1)}{4p(p-1)},
\end{align*}
where $k>0,1+\sqrt{2}<p<\frac{3+\sqrt{5}}{2}$.

\underline{Case 4.3.2.} $s\neq0$, $T$ is a rational triangle.

Take $s=2km,t=k(m^2-1),k>0,m>1$, then
\[w_1=k(m^2+1).\]
From the second equation of Eq. (\ref{4-1}), we have
\begin{align*}
r=\frac{2k((m^2+1)c-2m)}{c^2-1},~~w_2=\frac{k((m^2+1)c^2-4mc+m^2+1)}{c^2-1},
\end{align*}
where $c$ is a rational number.

From the third equation of Eq. (\ref{4-1}), we have
\[\left(\frac{km}{a}+bk(m^2-1)\right)^2+a^2k^2(m^2-1)^2=w_3^2.\]
Let
\begin{align*}
\frac{km}{a}+bk(m^2-1)=d(p^2-1),~~ak(m^2-1)=2dp,
\end{align*}
then
\begin{align*}
a=\frac{2dp}{k(m^2-1)},~~b=\frac{2p(p^2-1)d^2-k^2m(m^2-1)}{2dkp(m^2-1)},~~w_3=d(p^2+1).
\end{align*}
Eq. (\ref{4-1}) now becomes
\begin{equation}
\begin{split}
&2dp((p^2+1)d-k(m^2+1))c^2\\
&-k(m-1)(2p(m-1)d-k(m+1)(m^2+1))c\\
&-2p(p^2+1)d^2+4dkmp-2k^2m(m^2-1)=0,                                               \label{4-8}
\end{split}
\end{equation}
then
\[c=\frac{k(2(m-1)^2pd-(m^4-1)k)\pm\sqrt{\Delta(c)}}{4dp((p^2+1)d-k(m^2+1))},\]
where
\begin{align*}
\Delta(c)=&~(m^4-1)^2k^4-4dp(m^4-1)(m+1)^2k^3\\
&+4d^2p(4m(m^2-1)p^2+(m+1)^4p+4m(m^2-1))k^2\\
&-16d^3p^2(p^2+1)(m+1)^2k+16d^4p^2(p^2+1)^2.
\end{align*}
To get the rational solutions of Eq.\ (\ref{4-8}) with respect to $c$, it needs $\Delta(c)$ be a perfect square, say $w^2$. Since the discriminant $\Delta(c)$ is homogeneous with respect to $k$ and $d$, we can let $d=1$. Using a method described by Fermat (see \cite[p.~639]{Dickson}), it's easy to get the following solution
\begin{equation*}
\begin{split}
k=\frac{p^2+1}{m^2+1}
\end{split}
\end{equation*}
such that $\Delta(c)$ is a perfect square.

Let us consider the rational points on the quartic curve
\[\mathcal{C}: w^2=\Delta(c).\]
By the map $\varphi$:
\begin{align*}
X=&\frac{12pX(k,w)}{k^2},\\
Y=&\frac{432p^2(p^2+1)Y(k,w)}{k^3},
\end{align*}
where
\begin{align*}
X(k,w)=&~(4m(m^2-1)p^2+(m+1)^4p+4m(m^2-1))k^2\\
&-12p(p^2+1)(m+1)^2k+24p(p^2+1)^2-6(p^2+1)w,\\
Y(k,w)=&~(m^4-1)(m+1)^2k^3-2(4m(m^2-1)p^2+(m+1)^4p+4m(m^2-1))k^2\\
&+12p(p^2+1)(m+1)^2k-16p(p^2+1)^2-((m+1)^2k-4(p^2+1))w.
\end{align*}
We can transform $\mathcal{C}$ into the elliptic curve
\begin{align*}
\mathcal{E}_{(m,p)}:~Y^2=~X^3-432p^2A_4(m,p)X-3456p^3A_6(m,p),
\end{align*}
where
\begin{align*}
A_4(m,p)=&~4(m^2+3)(3m^2+1)(m^2-1)^2p^4-4(m-1)(3m^2-2m+3)(m+1)^5p^3\\
&+(25m^6-42m^5+119m^4-140m^3+119m^2-42m+25)(m+1)^2p^2\\
&-4(m-1)(3m^2-2m+3)(m+1)^5p+4(m^2+3)(3m^2+1)(m^2-1)^2,\\
A_6(m,p)=&~(2(m^2-1)(3m^2-2m+3)p^2-(m+1)^4p\\
&+2(m^2-1)(3m^2-2m+3))(8m(3m^2+2m+3)(m^2-1)^2p^4\\
&-4(m-1)(3m^2-2m+3)(m+1)^5p^3\\
&+(m^6+54m^5-49m^4+52m^3-49m^2+54m+1)(m+1)^2p^2\\
&-4(m-1)(3m^2-2m+3)(m+1)^5p+8m(3m^2+2m+3)(m^2-1)^2).
\end{align*}
By some calculations, we can check that the discriminant $\Delta(m,p)$ of $\mathcal{E}_{(m,p)}$ is not $0$, then $\mathcal{E}_{(m,p)}$ is non-singular.

It is easy to check that the elliptic curve $\mathcal{E}_{(m,p)}$ contains two rational points
\begin{align*}
&Q=(X_1(m,p),0),\\
&P=(X_2(m,p),864p^2(p^2+1)(m^4-1)(m+1)^2),
\end{align*}
where
\begin{align*}
X_1(m,p)=&-12p(2(m^2-1)(3m^2-2m+3)p^2\\
&-(m+1)^4p+2(m^2-1)(3m^2-2m+3)),\\
X_2(m,p)=&~12p(2(m^2-1)(3m^2+2m+3)p^2\\
&+(m+1)^4p+2(m^2-1)(3m^2+2m+3)).
\end{align*}
By the inverse map $\varphi^{-1}$, we have
\begin{align*}
k=-\frac{24p(p^2+1)k(X,Y)}{(X-X_1(m,p))(X-X_2(m,p))},
\end{align*}
where
\begin{align*}
k(X,Y)=&~Y+6p(m+1)^2(X+12p(2(m^2-1)(3m^2-2m+3)p^2\\
&-(m+1)^4p+2(m^2-1)(3m^2-2m+3))).
\end{align*}

By the group law, we get
\begin{align*}
[2]P=(&12(3(m^2-1)^2p^4-2(m^2-1)(3m^2-2m+3)p^3\\
&+(m+1)^2(7m^2-10m+7)p^2-2(m^2-1)(3m^2-2m+3)p+3(m^2-1)^2),\\
&-216(p^2+1)^2(m^2-1)^2((m^2-1)p^2-(3m^2-2m+3)p+m^2-1)).
\end{align*}
From the rational point $[2]P$, we have
\begin{align*}
k=\frac{4p}{m^2-1},
\end{align*}
then
\[c=-\frac{(m^2-1)p^2+8mp+m^2-1}{(m^2-1)p^2-4(m^2+1)p+m^2-1}.\]
By some calculations, we obtain
\begin{align*}
&r=-\frac{(p^2+1)((m^2-1)p^2-(4(m^2+1))p+m^2-1)}{(m^2-1)p^2-2(m-1)^2p+m^2-1},~~s=\frac{8mp}{m^2-1},~~t=4p,\\
&a=\frac{1}{2},~~b=\frac{(m^2-1)p^2-8mp-m^2+1}{4p(m^2-1)},
\end{align*}
where
\begin{align*}
1<p<2+\sqrt{3},~m>1~~\text{or}~~p>2+\sqrt{3},~1<m<\frac{\sqrt{(p^2-4p+1)(p^2+4p+1)}}{p^2-4p+1},
\end{align*}
this is a parametric solution of Eq. (\ref{4-1}).

\section{areas and perimeters are proportional}
In this section, we consider the problem that rational triangle pairs or rational triangle-parallelogram pairs with areas and perimeters in fixed proportions $(\alpha,\beta)$ respectively, where $\alpha$ and $\beta$ are positive rational numbers. The existed results can be found in \cite{Li-Zhang1,Li-Zhang2,Li-Zhang3}. We only give the Diophantine systems corresponding to Eq. (\ref{2}) and Eq. (\ref{4-1}), respectively.

\underline{Case 5.1.} Rational triangle pairs.

The transformation matrix $M_f$ can be taken as
$$
M_f=\left[\begin{matrix}
\frac{\alpha}{a} & b \\
0 & a \\
\end{matrix}\right],$$
where $a\in\mathbb{Q^*},b\in\mathbb{Q}$. The corresponding Diophantine system is
\begin{equation}
\begin{cases}
\begin{aligned}
&s^2+t^2=w_1^2,\\
&\left(s-r\right)^2+t^2=w_2^2,\\
&\left(\frac{\alpha s}{a}+bt\right)^2+\left(at\right)^2=w_3^2,\\
&\left(\frac{\alpha(s-r)}{a}+bt\right)^2+(at)^2=w_4^2,\\
&\frac{\alpha r}{a}+w_3+w_4=\beta\left(r+w_1+w_2\right),                                  \label{5-1}
\end{aligned}
\end{cases}
\end{equation}
where $b\in\mathbb{Q}$, $s\in\mathbb{Q}_{\geq0}$, $\alpha,\beta,a,r,t,w_1,w_2,w_3,w_4\in\mathbb{Q^*}$.

\underline{Case 5.2.} Rational triangle and parallelogram pairs.

The transformation matrix $M_f$ can be taken as
$$
M_f=\left[\begin{matrix}
\frac{\alpha}{2a} & b \\
0 & a \\
\end{matrix}\right],$$
where $a\in\mathbb{Q^*},b\in\mathbb{Q}$. The corresponding Diophantine system is
\begin{equation}
\begin{cases}
\begin{aligned}
&s^2+t^2=w_1^2,\\
&\left(s-r\right)^2+t^2=w_2^2,\\
&\left(\frac{\alpha s}{2a}+bt\right)^2+\left(at\right)^2=w_3^2,\\
&2\left(\frac{\alpha r}{2a}+w_3\right)=\beta(r+w_1+w_2),                                               \label{5-2}
\end{aligned}
\end{cases}
\end{equation}
where $b\in\mathbb{Q}$, $s\in\mathbb{Q}_{\geq0}$, $\alpha,\beta,a,r,t,w_1,w_2,w_3\in\mathbb{Q^*}$.

The solutions of Eq. (\ref{5-1}) and Eq. (\ref{5-2}) can be discussed similarly to Eq. (\ref{2}) and Eq. (\ref{4-1}), respectively.

     \end{document}